\documentclass[11pt, a4paper]{article}
\def\dj{d\kern-0.4em\char"16\kern-0.1em}

\usepackage{a4wide}

\usepackage{amssymb,amsmath,amsthm}
\usepackage{graphicx,enumerate}
\usepackage{ragged2e}
\usepackage{xcolor}
\usepackage[utf8]{inputenc}
\usepackage[T1]{fontenc}
\usepackage[english]{babel}
\usepackage{subcaption}
\usepackage{amssymb, amsmath, latexsym, amsthm, enumerate, amsfonts}
\usepackage[labelsep=period]{caption}

\usepackage[font=small]{caption}
\usepackage{mwe}

\usepackage[pdftex, pdfstartview=FitH]{hyperref}

\numberwithin{equation}{section}

 \newtheorem{thm}{Theorem}[section]
 \newtheorem{cor}[thm]{Corollary}
 \newtheorem{con}[thm]{Conjecture}
  
 \newtheorem{lem}[thm]{Lemma}

 \newtheorem{prop}[thm]{Proposition}
 \theoremstyle{definition}
 
 \theoremstyle{remark}

\newcommand\cF{{\mathcal F}}
\newcommand\cG{{\mathcal G}}
\newcommand\cH{{\mathcal H}}
\newcommand\cP{{\mathcal P}}
\newcommand\cU{{\mathcal U}}

\newcommand{\La}{\mathop{}\!\mathrm{La}}

\newcommand{\vex}{\mathop{}\!\mathrm{vex}}

\begin{document}

\title{A note on vertex Turán problems in the Kneser cube}
\author{D\'aniel Gerbner\thanks{HUN-REN Alfr\'ed R\'enyi Institute of Mathematics, Budapest, 1053, Re\'altanoda utca 13-15, e-mail: \href{mailto:gerbner@renyi.hu}{gerbner@renyi.hu}, Research supported by NKFIH under grants FK132060 and KKP-133819.} \and Bal\'azs Patk\'os\thanks{HUN-REN Alfr\'ed R\'enyi Institute of Mathematics, Budapest, 1053, Re\'altanoda utca 13-15, e-mail:  \href{mailto:patkos@renyi.hu}{patkos@renyi.hu} Research supported by NKFIH under grant FK132060.}}
\date{}

\maketitle

\begin{abstract}
    The Kneser cube $Kn_n$ has vertex set $2^{[n]}$ and two vertices $F,F'$ are joined by an edge if and only if $F\cap F'=\emptyset$. For a fixed graph $G$, we are interested in the most number $\vex(n,G)$ of vertices of $Kn_n$ that span a $G$-free subgraph in $Kn_n$. We show that the asymptotics of $\vex(n,G)$ is $(1+o(1))2^{n-1}$ for bipartite $G$ and $(1-o(1))2^n$ for graphs with chromatic number at least 3. We also obtain results on the order of magnitude of  $2^{n-1}-\vex(n,G)$ and $2^n-\vex(n,G)$ in these two cases. In the case of bipartite $G$, we relate this problem to instances of the forbidden subposet problem.
\end{abstract}

\section{Introduction}

Graph Tur\'an problems are among the most studied problems in extremal graph theory \cite{FS}. They ask for how large a subgraph $F$ of a given host graph $H$ can be if $F$ does not contain $G$ as a subgraph, i.e. the input of the problem is ($G$,$H$). Largeness is usually measured by the number of edges, but there are other parameters like spectral radius. Recently there have been interest \cite{AKSZ,Getal,JT,LH} in \textit{vertex Tur\'an problems}, where one is interested in the most number of vertices that a vertex subset $U$ of $H$ can have such that the induced subgraph $H[U]$ is $G$-free. 

One particular host graph where this problem has been studied is the \textit{Kneser graph} $Kn(n,k)$ with vertex set $\binom{[n]}{k}:=\{F\subseteq [n]:\, |F|=k\}$, where $[n]=\{1,2,\dots,n\}$, and $F,F'$ are connected by an edge if and only if $F$ and $F'$ are disjoint. Vertex Tur\'an problems in the Kneser graph are usually phrased as extremal finite set theory problems. The first such result is the Erd\H os-Ko-Rado theorem \cite{EKR}, which determines the largest family of pairwise intersecting $k$-sets. In other words, the Erd\H os-Ko-Rado theorem determines the largest independent set, i.e., the largest $K_2$-free vertex subset of $Kn(n,k)$. The Erd\H os matching problem \cite{E2} deals with the largest family without $r+1$ pairwise disjoint sets, i.e., the largest set of vertices without $K_{r+1}$. For results on this problem, see \cite{FK2.5} and the citations in it. Gerbner, Lemons, Palmer, Patk\'os and Sz\'ecsi \cite{glppsz} studied the largest $\ell$-almost intersecting families, where each set is disjoint with at most $\ell$ other sets, i.e., the largest $S_{\ell+1}$-free sets of vertices. Katona and Nagy \cite{KN} studied $(s,t)$-union intersecting families, which have the property that the union of any $s$ sets intersect the union of any $t$ sets, i.e., $K_{s,t}$-free sets. Alishahi and Taherkhani \cite{AT} united these notions and determined for every $F$ the largest $F$-free set of vertices in the Kneser graph $Kn(n,k)$, provided $n$ is sufficiently large. More results in this topic can be found in \cite{Getal2,T}.

In each of the above specific problems on set families, the non-uniform version has also been studied, but not for the general problem. Here we address this.
Our host graph is the Kneser cube $Kn_n$ with vertex set $2^{[n]}:=\{F\subseteq [n]\}$, where $F,F'$ are connected by an edge if and only if $F$ and $F'$ are disjoint. We denote by $\vex(n,G)$ the most number of vertices of $Kn_n$ that induce a $G$-free subgraph in $Kn_n$. Note that $\vex(n,G)\ge 2^{n-1}$ if $G$ contains an edge, as shown by the family of sets each containing a fixed element. Therefore, the order of magnitude of $\vex(n,G)$ is known, and the more interesting question may be to determine or bound the order of magnitude of $2^n-\vex(n,G)$ or $\vex(n,G)-2^{n-1}$.

If $G=K_{r+1}$, then the problem is equivalent to the non-uniform Erd\H os matching problem. Kleitman \cite{K} obtained a general and very nice upper bound on $\vex(n, K_{r+1})$. Sharpenings for special cases of $r$ and the residue class of $n$ modulo $r+1$ have been obtained in (among others) \cite{FK1,FK2,FK3}. A simple lower bound is obtained by taking each set of size larger than $n/(r+1)$; all the above results show upper bounds close to this.

Gerbner, Lemons, Palmer, Patk\'os and Sz\'ecsi \cite{glppsz} determined $\vex(n,S_2)$ and $\vex(n,S_3)$ exactly.
Katona and Nagy \cite{KN} also determined $\vex(n,C_4)$ exactly and for $t\ge 4$ they showed that $\vex(n,S_t)-2^{n-1}=(\frac{1}{2}+o(1))\binom{n}{\lfloor n/2\rfloor}$. For $t\ge s\ge 2$, Katona and Nagy showed that $\vex(n,K_{s,t})-2^{n-1}=(1+o(1))\binom{n}{\lfloor n/2\rfloor}$.

We shall consider all graphs $G$ and the corresponding vertex Tur\'an problem in the Kneser cube. It turns out that there are three major graph classes with respect to the value of $\vex(n,G)$: matchings, other bipartite graphs and graphs with chromatic number at least 3. The following proposition about matchings is rather a simple exercise. We include its proof for completeness in Section \ref{secbip}, the interested reader is welcome to figure out the proof on its own. $M_k$ denotes the matching with $k$ edges.

\begin{prop}\label{match}
$\vex(n,M_{k+1})=2^{n-1}+k$.    
\end{prop}

The situation is somewhat different for bipartite graphs other than matchings. The asymptotics of $\vex(n,G)$ is still $2^{n-1}$, but the second order term is different. Our next result determines the order of magnitude of $2^{n-1}-\vex(n,G)$.

\begin{thm}\label{bip}
    For any bipartite graph $G$ that is not a matching, we have $$2^{n-1}+\left(\frac{1}{2}-o(1)\right)\binom{n}{n/2}\le\vex(n,G)\le 2^{n-1}+(1+o(1))\binom{n}{\lfloor n/2\rfloor}.$$
\end{thm}

Both bounds of Theorem \ref{bip} can be tight. In Section \ref{secbip}, we shall relate the problem of finding $\vex(n,G)$ for bipartite graphs $G$ to a forbidden subposet problem (see the definitions there). Whenever the corresponding problem is solved, we will be able to determine $\vex(n,G)$ up to a $o(\binom{n}{n/2})$ term and whenever that forbidden subposet problem is settled, the coefficient of $\binom{n}{n/2}$ is always either $1/2$ or 1. 

\vskip 0.2truecm

Next, we shall turn our attention to non-bipartite graphs. We introduce the following notation: $\binom{n}{\le m}:=\sum_{i=0}^{m}\binom{n}{i}$ and $\binom{n}{\ge m}:=\sum_{i=m}^{n}\binom{n}{i}$. The following simple proposition shows that if the chromatic number of $G$ is at least 3, then almost all vertices can be selected to obtain a $G$-free induced subgraph of $Kn_n$. 

\begin{prop}\label{constr}
    For any graph $G$ with chromatic number at least 3, let $\ell=\ell(G)$ be the length of the shortest odd cycle of $G$ with $\ell=2k+1$. Then $2^n-\vex(n,G)\le \binom{n}{\le \lfloor\frac{kn}{2k+1}\rfloor}$.
\end{prop}

Our final theorem establishes that for odd cycles, the correct order of magnitude of $2^n-\vex(n,C_{2k+1})$ is given by the bound of Proposition \ref{constr}. In the next theorem and later, $\log$ denotes the natural logarithm.

\begin{thm}\label{cycle}
For any positive integer $k$, we have $\vex(n,C_{2k+1})\le \binom{n}{\ge \lceil\frac{kn}{2k+1}\rceil-\lceil2(k+1)\log(2k)\rceil}$. Therefore, $2^n-\vex(n,C_{2k+1})=\Theta_k(\binom{n}{\le \lfloor\frac{kn}{2k+1}\rfloor})$.
\end{thm}

\section{Proofs}
\subsection{Bipartite graphs}\label{secbip}

We start with the proof of Proposition \ref{match}.

\begin{proof}[Proof of Proposition \ref{match}]
    Let $\cF\subseteq 2^{[n]}$ be such that $Kn_n[\cF]$ is $M_{k+1}$-free. To avoid a copy of $M_{k+1}$, out of every pair of complement sets at most $k$ can exist with both sets belonging to $\cF$. This shows the upper bound. And for any family $\cF$ containing exactly one of $F,[n]\setminus F$ for $2^{n-1}-k$ pairs of complements and both of them for the remaining $k$ pairs, the induced subgraph $Kn_n[\cF]$ is $M_k$ with isolated vertices. 
\end{proof}

To prove Theorem \ref{bip}, we need some auxiliary definitions and results. We say that a family $\cG$ is a copy of a poset $P=(P,\prec)$ if there exists a bijection $b:P\rightarrow \cG$ such that $p\prec q$ implies $b(p)\subseteq b(q)$ for any $p,q\in P$. A family $\cF$ is $P$-free if it does not contain any copy of $P$. The problem of finding the size of the largest $P$-free family $\cF\subseteq 2^{[n]}$, denoted by $\La(n,P)$, is called the forbidden subposet problem (for a survey see \cite{GL} and Chapter 7 of \cite{GP}). If set $\cP$ of posets is forbidden, then $\La(n,\cP)$ denotes the maximum size of a family $\cF\subseteq 2^{[n]}$ that is $P$-free for all $P\in \cP$. It follows from an old result of Erd\H os \cite{E} that for any poset $P$, we have $\La(n,P)\le (|P|-1)\cdot \binom{n}{\lfloor n/2\rfloor}$, but there exist stronger results (and used to be a stronger conjecture - see at the end of this subsection). We will use a special case of the following theorem of Bukh \cite{B}. The \textit{Hasse diagram} of a poset $(P,\prec)$ is the oriented graph with vertex set $P$ and $\overrightarrow{pq}$ is an arc if $p\prec q$ and there is no $z\in P$ with $p\prec z\prec q$. A poset is called a \textit{tree poset} if its Hasse diagram is a tree. The height $h(P)$ of a poset $P$ is the size of its longest chain.

\begin{thm}[Bukh \cite{B}]\label{bukh}
    For any tree poset $T$, we have $\La(n,T)=(h(T)+o(1))\binom{n}{n/2}$.
\end{thm}

Let $G=(A,B,E)$ be a bipartite graph. We define $P_{G,A}$ to be the poset with Hasse diagram $G$ and each edge oriented towards its end-vertex in $A$. The poset $P_{G,B}$ is the \textit{dual poset} of $P_{G,A}$, that is, the poset with the same Hasse diagram but with all edges oriented in the opposite direction.

\begin{lem}\label{equiv}
Let $G$ be a connected bipartite graph, and let $\cF\subseteq 2^{[n]}$ be such that $Kn_n[\cF]$ is $G$-free. Then $\cF_{sym}:=\cF\cap \cF^c\cap \bigcup_{i=-n^{2/3}}^{n^{2/3}}\binom{[n]}{\lfloor n/2\rfloor+i}$ is $P_{G,A}$-free if $n$ is large enough.
\end{lem}

\begin{proof}
    Suppose towards a contradiction that $\cF_{sym}$ contains a copy $\cH$ of $P_{G,A}$.  Let $\cH_A$ and $\cH_B$ be the parts of $\cH$ that correspond to vertices of $A$ and $B$, respectively. First we show that we cannot have $F\in \cH_a$ and $F^c\in\cH_B$. This is because $G$ is connected, and therefore $\cH$ is connected. As a result, there is a path $v_1v_2\dots v_{2k}$ in $\cH$ with $v_1=F$, $v_{2k}=F^c$ and $2k\le |V(G)|$. Observe that the Hamming distance between $v_i$ and the complement of $v_{i+1}$ is at most $2n^{2/3}$ for any $i$. Therefore, the Hamming distance between $v_i$ and $v_{i+2}$ is at most $4n^{2/3}$. This implies that the Hamming distance between $v_1$ and $v_{2k-1}$ is at most $2(2k-2)n^{2/3}$, thus $v_{2k-1}$ is not disjoint from the complement of $v_1$, a contradiction.  

    Consider the family $\cH_B\cup \cH^c_A$. If there is a relation, i.e., a directed edge between $H\in \cH_B$ and $H'$ in $G_{P,A}$, then the edge goes to $H'$, i.e., $H\subset H'$. This implies $H\cap H'^c=\emptyset$, and thus $H$ and $H'^c$ are connected by an edge. So by definition of $P_{G,A}$, the family $\cH_B\cup \cH^c_A$ induces a copy of $G$ in $Kn_n$, a contradiction completing the proof. 
\end{proof}

Now we are ready to prove Theorem \ref{bip}.

\begin{proof}[Proof of Theorem \ref{bip}]
    Let $G$ be any bipartite graph that has at least one vertex of degree at least 2, i.e. not a matching. First we show the lower bound $\vex(n,G)\ge 2^{n-1}+\frac{1}{2}\binom{n}{n/2}$ if $n$ is even, and $\vex(n,G)\ge 2^{n-1}+\binom{n-1}{\lfloor n/2\rfloor-1}$ if $n$ is odd. In the former case the family of all the sets of size at least $n/2$ induce a matching in $Kn_n$, and the same holds in the odd $n$ case for the family consisting of all the sets of size larger than $n/2$ together with the sets of size $\lfloor n/2\rfloor$ sharing a common element.

    To see the upper bound, let $\cF\subseteq 2^{[n]}$ be a family with $Kn_n[\cF]$ being $G$-free. Observe first that by the well-known Chernoff bound, $|\binom{[n]}{\le n/2-n^{2/3}}\cup \binom{[n]}{\ge n/2+n^{2/3}}|=o(\binom{n}{n/2})$. We claim that $\cF_{sym}:=\cF\cap \cF^c\cap \bigcup_{i=n/2-n^{2/3}}^{n/2+n^{2/3}}\binom{[n]}{i}$ have size at most $\La(n,P_{G,A})+2|G|$. Let $C^1,C^2,\dots, C^m$ be the components of $G$ and $A^j\subset A$ is the part of $C^j$ that belongs to $A$. If $|\cF_{sym}|>\La(n,P_{G,A})+2|G|$, then by  Lemma \ref{equiv} and by $\La(n,P_{C^j,A^j})\le \La(n,P_{G,A})$ for any $1\le j\le m$, we obtain a copy of $C^j$. Removing this copy and the complement sets, the remaining family is still large enough to find the second component, and so on. As each time we throw away $2|C^j|$ sets, therefore even before finding the last component we still have more than $|\cF|-2|G|>\La(n,P_{G,A})$ sets. 
    So we found a copy of $G$ and this contradiction yields $|\cF|\le 2^{n-1}+o(\binom{n}{n/2}+\frac{1}{2}\La(n,P_{G,A})+2|G|$.

Observe that if $G$ is bipartite, then $P_{G,A}$ has height 2 (all chains have length at most 2) and all such posets are subposets of $K_{s,1,t}$ for appropriately chosen $s$ and $t$, where $K_{s,1,t}$ is the so-called complete 3-level poset with elements $u_1,u_2,\dots,u_s,v,w_1,w_2,\dots,w_t$ and $u_i\prec v\prec w_j$ for any $1\le i\le s,1\le j\le t$ being all its relations. As $K_{s,1,t}$ is a tree poset, Theorem \ref{bukh} implies that for all bipartite $G$, we have $\La(n,P_{G,A})\le (2+o(1))\binom{n}{n/2}$. Plugging this in the upper bound at the end of the previous paragraph, we obtain the statement of the theorem.
\end{proof}

In obtaining the upper bound of Theorem \ref{bip}, we used the inequality $\vex(n,G)\le 2^{n-1}+\frac{1}{2}\La(n,P_{G,A})+o(\binom{n}{n/2})$. In a similar fashion, $\vex(n,G)\le 2^{n-1}+\frac{1}{2}\La(n,P_{G,B})+o(\binom{n}{n/2})$ and $$\vex(n,G)\le 2^{n-1}+\frac{1}{2}\La(n,\{P_{G,A},P_{G,B}\})+o(\binom{n}{n/2})$$ hold, where $\La(n,\{P_1,P_2\})$ denotes the maximum size of a family $\cF\subseteq 2^{[n]}$ that is both $P_{G,A}$- and $P_{G,B}$-free. Also, if $\La_{sym}(n,P)$ denotes the maximum size of a $P$-free family $\cF\subseteq 2^{[n]}$ that is complement closed, i.e. $\cF^c=\cF$, then the above reasoning yields
\[
\vex(n,G)\le 2^{n-1}+\frac{1}{2}\La_{sym}(n,P_{G,A})+o(\binom{n}{n/2}).
\]
We conjecture that this upper bound gives the correct asymptotics of $\vex(n,G)$ for all bipartite graphs that are not matchings.

\begin{con}
    For any bipartite graph $G$ that is not a matching, we have $\vex(n,G)=2^{n-1}+\frac{1}{2}\La_{sym}(n,P_{G,A})+o(\binom{n}{n/2})$.
\end{con}

 There is a very important poset parameter that is relevant in forbidden subposet problems. For a poset $P$, let $e(P)$ denotes the largest integer $k$ such that for any $n$ and $j<n$ the family $\bigcup_{i=1}^k\binom{[n]}{j+i}$ is $P$-free. In words, $e(P)$ denote the maximum number of consecutive levels of the Boolean cube that we can have without creating a copy of $P$. This means that if we consider the family containing the middle $e(P)$ levels of the Boolean cube, then we obtain the lower bound $(e(P)+o(1))\binom{n}{n/2}\le \La (n,P)$.
There used to be a very important conjecture of Bukh \cite{B} and Griggs, Lu \cite{GL2} in the area of forbidden subposet problems. It stated that this construction is asymptotically optimal for any poset $P$, and so $\La(n,P)=(e(P)+o(1))\binom{n}{n/2}$. This has been very recently disproved by Ellis, Ivan, and Leader \cite{EIL}, but their smallest counterexample is the Boolean poset $B_4$ of dimension 4 with $e(P)=4$. Our posets of interest all have height 2.

What bearings does the parameter $e(P)$ have on $\vex(n,G)$ for bipartite graphs $G$? As we discussed in the proof of Theorem \ref{bip}, $P_{G,A}$ has height 2 and thus it is a subposet of $K_{s,1,t}$ for $s,t$ large enough, and thus $e(P_{G,A})\le 2$, and so $e(P_{G,A})$ is either 1 or 2. Note that $e(P_{G,A})=e(P_{G,B})$ by the definition of $e$ and $(e(P_{G,A})+o(1))\binom{n}{n/2}\le \La(n,\{P_{G,A},P_{G,B}\})$. Can we turn the middle level construction showing this lower bound into a construction that gives the analogous lower bound for $\vex(n,G)$? The answer is positive, but let us start with the simpler case $e(G_{P,A})=1$. Putting together the lower bound of Theorem \ref{bip} and the upper bound $\vex(n,G)\le 2^{n-1}+o(\binom{n}{n/2})+\frac{1}{2}\La(n,P_{G,A})$, we obtain the following corollary for those bipartite graphs $G$ for which the generally false conjecture by Bukh and Griggs-Lu holds for $P_{G,A}$ with $e(P_{G,A})=1$.

\begin{cor}
    If for a bipartite graph $G=(A,B,E)$ that is not a matching, we have $e(P_{G,A})=1$ and $\La(n,P_{G,A})=(1+o(1))\binom{n}{n/2}$, then $\vex(n,G)=2^{n-1}+(1/2+o(1))\binom{n}{n/2}$ holds.
\end{cor}

Let us now consider bipartite graphs with $e(P_{G,A})=2$. This means that the union of the two middle levels does not contain $P_{G,A}$ nor $P_{G,B}$. For a a family $\cF\subseteq 2^{[n]}$ let us write $\cU(\cF):=\{G\subseteq [n]: \exists F\in \cF, F\subseteq G\}$. Let us introduce the following families:
\[\cF^{even}_0=\{F\subset [n]:1\in F, |F|=n/2-1,n/2\}\cup \{F\subset [n]:1\notin F, |F|=n/2,n/2+1\},
\]
\[
\cF^{odd}_0=\binom{[n]}{\lfloor n/2\rfloor}\cup \binom{n}{\lceil n/2\rceil}, \hskip 0.5truecm \cF^{odd}=\cU(\cF^{odd}_0),  \hskip 0.5truecm \cF^{even}=\cU(\cF^{even}_0).
\]
Now $\cF^{odd}_0$ cannot contain $P_{G,A}$ nor $P_{G,B}$ by definition of $e$, while $\cF_0^{even}$ has this property as there is no containment between $\{F\subset [n]:1\in F, |F|=n/2-1,n/2\}$ and $\{F\subset [n]:1\notin F, |F|=n/2,n/2+1\}$, and separately these two subfamilies are two consecutive levels of Boolean cubes of dimension $n-1$. Next, we claim that $Kn_n[\cF^{odd}]$ and $Kn_n[\cF^{even}]$ are $G$-free. Indeed, vertices corresponding to sets in $\cF^{even}\setminus \cF^{even}_0$ and $\cF^{odd}\setminus \cF^{odd}_0$ are isolated in $Kn_n[\cF^{even}]$ and $Kn_n[\cF^{odd}]$, respectively. So a copy of $G$ should be in $Kn_n[\cF^{even}_0]$ or in $Kn_n[\cF^{odd}_0]$. Note that if $G$ has components $C^1,C^2,\dots,C^r$, then $e(P_{G,A})=\max\{e(P_{C^i,A_i}):1\le i\le r\}$, so there exists a $C^i$ with $e(P_{C^i,A_i})=2$. A copy of $C^i$ in $Kn_n[\cF^{even}]$ or $Kn_n[\cF^{odd}]$ would yield a copy of $P_{C^i,A_i}$ in $Kn_n[\cF^{even}]$ or $Kn_n[\cF^{odd}]$ as in the proof of Theorem \ref{bip}. But that contradicts $e(P_{C^i,A_i})=2$. So we obtained a lower bound $\vex(n,G)\ge 2^{n-1}+(1+o(1))\binom{n}{n/2}$ for arbitrary bipartite graphs $G$ with $e(P_{G,A})=2$. Together with the upper bound of Theorem \ref{bip}, it yields the following corollary.

\begin{cor}\label{e2}
    If for a bipartite graph $G=(A,B,E)$, we have $e(P_{G,A})=2$, then $\vex(n,G)=2^{n-1}+(1+o(1))\binom{n}{n/2}$ holds.
\end{cor}

Let us finish this subsection with two remarks on the smallest/simplest examples of bipartite graphs $G$ for which $\vex(n,G)$ is not known up the constant term of $\binom{n}{n/2}$. Among even cycles we have $e(P_{C_4,A})=2$, and for any $t\ge 3$, we have $e(P_{C_{2t},A})=1$. The corresponding posets have special names: $P_{C_4,A}$ is the butterfly poset $\bowtie$ and $P_{C_{2t},A}$ for $\ge 3$ is the crown poset also denoted by $C_{2t}$. (Note that $P_{C_{2t},A}=P_{C_{2t},B}$ for any $t\ge 2$.)  Corollary \ref{e2} gives the asymptotics of $\vex(n,C_4)$, but, as mentioned in the introduction, its \textit{exact value} was determined by Katona and Nagy \cite{KN}. The equality $\La(n,P)=(e(P)+o(1))\binom{n}{n/2}$ is known to be true for $C_{2t}$ (for $t$ even proved by Griggs and Lu \cite{GL2}, for $t$ odd proved by Lu \cite{L}),
with the exception of $C_6$ and $C_{10}$. So these two graphs are natural targets of future research.

As we have mentioned several times, Ellis, Ivan, and Leader \cite{EIL} disproved the conjecture $\La(n,P)=(e(P)+o(1))\binom{n}{n/2}$. They constructed a family $\cF\subseteq \binom{[n]}{n/2}$ of size $0.29\binom{n}{n/2}$ that does not contain all six sets of any 'middle 4-cube', i.e. for any $S\subset  [n]$ with $|S|=n/2-2$, $|T|=n/2+2$, we have $\{G: S\subset G\subset T, |G|=n/2\}\not\subset \cF$. Adding sets of size $n/2\pm 1,2$ shows $\La(n,B_4)\ge (4.29+o(1))\binom{n}{n/2}$, where $B_4$ is the Boolean lattice of dimension 4. Also, let $G_4$ be the bipartite incidence graph of the complete graph $K_4$, that is $G_4=(A,B,E)$ with $A=E(K_4),B=V(K_4)$ and $(e,v)\in E$ if and only if $v\in e$. Then the above $\cF$ together with all the sets of size $n/2-1$ shows that $\La(n,P_{G_4,A})\ge (1.29+o(1))\binom{n}{n/2}$. Similarly, $\cF\cup \binom{[n]}{n/2+1}$ shows $\La(n,P_{G_4,B})\ge (1.29+o(1))\binom{n}{n/2}$. On the other hand, the first family contains a copy of $P_{G_4,B}$ while the second contains a copy of $P_{G_4,A}$. If the original conjecture would have been true, then the parameter $\La(n,P)$ would have been \textit{principal}, i.e. $\La(n,\mathcal{P})=(1+o(1))\min \{\La(n,P): P\in \mathcal{P}\}$ for any family $\mathcal{P}$ of posets. In view of the above, we conjecture that this is not true, in particular we make the following conjecture.

\begin{con}
    $\La(n,P_{G_4,A})-\La(n,\{P_{G_4,A},P_{G_4,B}\})=\Omega(\binom{n}{n/2})$.
\end{con}

Let us point out that by a result of Griggs and Lu (Theorem 1.4 in \cite{GL2}) $\La(n,P_{G_4,A})\le (1+\sqrt{2/3}+o(1))\binom{n}{n/2}$ and so if the limit $\lim_n\frac{\La(n,P_{G_4,A})}{\binom{n}{n/2}}$ exists, then its value is a non-integer.

\subsection{Non-bipartite graphs}

To address the problem of estimating $\vex(n,G)$ for non-bipartite graphs $G$, we start with the proof of Proposition \ref{constr}.

\begin{proof}[Proof of Proposition \ref{constr}]
    Consider the family $\cF_k=\{F\subseteq [n]: |F|>\frac{k}{2k+1}n\}$. We claim that $Kn_n[\cF_k]$ is $C_{2k+1}$-free. Indeed, for any $F_1,F_2,\dots, F_{2k+1}\in \cF_k$, we have $\frac{\sum_{i=1}^{2k+1}|F_i|}{n}>(2k+1)\cdot \frac{k}{2k+1}n=k$, and thus there exists an element $x\in [n]$ that belongs to at least $k+1$ $F_i$s. Therefore these sets form an independent set in $Kn_n$, which does not exist in $C_{2k+1}$. The size of $\cF_k$ is $2^n-\binom{n}{\le \frac{k}{2k+1}n}$. This completes the proof of the proposition. 
\end{proof}

Now we turn to the proof of Theorem \ref{cycle}. We will use Katona's cycle method \cite{kat}. Let $\pi$ be a cyclic permutation of $[n]$. We call a subset $S$ of $[n]$ an \textit{interval with respect to $\pi$} if $S=\{\pi(i),\pi(i+1),\pi(i+2)\dots,\pi(i+j)$ for some $1\le i\le n, 0\le j\le n-1$, where addition is modulo $n$, so for example $\{\pi(n-1),\pi(n),\pi(1)\}$ is an interval. We will omit "with respect to $\pi$" whenever it is clear from context. For a family $\cF\subseteq 2^{[n]}$ and a cyclic permutation $\pi$ of $[n]$, we define $\cF_\pi=\{F\in \cF: F ~\text{is an interval with respect to}\ \pi\}$. We shall also need the following weight function: $w_n:[n]\rightarrow \mathbb{N}$ with $w_n(i)=\binom{n}{i}$. The weight $w_n(F)$ of a set $F$ is $w_n(|F|)$, and we let $w_n(\cF)=\sum_{F\in \cF}w_n(F)$, and we will sometimes omit the subscript $n$ if it is clear from context. The following lemma reduces our problem "to the cycle".

\begin{lem}\label{katona}
    Let $\mathbf{G}$ be a set of pair $(\cG,\pi)$ with $\cG$ being a family of intervals with respect to $\pi$ and $\pi$ is a cyclic permutation of $[n]$. Assume that $w(\cG)\le n\cdot B$ for every $\cG$ for which there exists a $\pi$ such that $(\cG,\pi)\in \mathbf{G}$. Suppose further that for some $\cF\subseteq 2^{[n]}$ it holds that for any cyclic permutation $\pi$, we have $(\cF_\pi,\pi)\in \mathbf{G}$. Then $|\cF|\le B$.
\end{lem}

This lemma is well-known, we include a proof for the sake of completeness.

\begin{proof}
    Let us count $\sum_{\pi}w(\cF_\pi)$ in two ways, where the summation runs over all cyclic permutation $\pi$ of $[n]$. By the assumptions of the lemma, it is at most $(n-1)!\cdot nB=n!B$, as the number of cyclic permutations is $(n-1)!$. On the other hand, every $F\subseteq [n]$ is an interval with respect to exactly $|F|!(n-|F|)!$ cyclic permutations, so the sum equals $$\sum_{F\in \cF}|F|!(n-|F|)!\cdot w(F)=\sum_{F\in \cF}n!.$$
    We obtained
    \[
    |\cF|\cdot n!=\sum_{F\in \cF}n!\le n!B,
    \]
    and dividing by $n!$ yields the statement of the lemma.
\end{proof}

Let us continue with the proof of Theorem \ref{cycle}.

\begin{proof}[Proof of Theorem \ref{cycle}]
Based on Lemma \ref{katona}, it is enough to prove that for any family $\cG$ of intervals such that $Kn_n[\cG]$ is $C_{2k+1}$-free, we have $$w_n(\cG)\le n\sum_{i=\lceil\frac{kn}{2k+1}\rceil-\lceil2(k+1)\log(2k)\rceil}^n\binom{n}{i}.$$
For any $j\le \frac{kn}{2k+1}$ let $m(j)$ denote the maximum number $m$ such that $\lfloor \frac{kn}{2k+1}\rfloor-j+2km\le kn$. Observe that $m(j)$ is monotone increasing in $j$ and $m(j+2k)=m(j)+1$. Also, $m(j)\le n/2$ with equality only if $n$ is even and $j=\lfloor \frac{kn}{2k+1}\rfloor$. We also have $m(j)\ge \lceil\frac{kn}{2k+1}\rceil-\lceil(k+1)\log(2k)\rceil$.

Consider the following mapping: for any interval $F$ of size $\lfloor \frac{kn}{2k+1}\rfloor-j$, let $f(F)=\{F_1,F_2,\dots,F_{2k}\}$, where all $F_i$s are intervals of size $m(j)$ such that $F_1$ is the interval with left endpoint immediately after $F$'s right endpoint, $F_2$ is the interval with left endpoint immediately after $F_1$'s right endpoint, and so on. 

By definition and by the sizes of the intervals $F$ is disjoint with $F_1$, $F_i$ is disjoint with $F_{i+1}$ for all $1\le i \le 2k-1$. Also, by the definition of $m(j)$, $F$ is disjoint with $F_{2k}$. Indeed, if $F=\{1,\dots, \lfloor \frac{kn}{2k+1}\rfloor-j\}$, then the last element of $F_{2k}$ is $\lfloor \frac{kn}{2k+1}\rfloor-j+2km\le kn$, and the first element is more than $\lfloor \frac{kn}{2k+1}\rfloor-j$ if $n$ is sufficiently large.
Therefore, $F$ and the sets in $f(F)$ form a $C_{2k+1}$ in $Kn_n$. This means that for every $G\in \cG$ at least one set of $f(G)$ does not belong to $\cG$.

Observe that if $j$ is fixed, then an interval $F$ of size $m(j)$ belongs to exactly $2k$ images $f(G)$ with $|G|=\lfloor \frac{kn}{2k+1}\rfloor -j$. As mentioned above, for every $m>\frac{kn}{2k+1}$ there are exactly $2k$ values of $j$ for which $m(j)=m$. Suppose there are $M$ intervals of size $m$ missing from $\cG$. Then for any $j$ with $m(j)=m$, $\cG$ can contain at most $2kM$ intervals of size $\lfloor \frac{kn}{2k+1}\rfloor -j$. Therefore, if
\[
2k\sum_{j:m(j)=m}w_n\left(\left\lfloor \frac{kn}{2k+1}\right\rfloor -j\right)\le w_n(m),
\]
then the weight of all the $n$ intervals of length $m$ is larger than the weight of all the $j$-element intervals with $m(j)=m$ plus the weight of the $m$-element intervals in $\cG$. In other words, to achieve largest possible weight, it is not worth putting any interval of length $\lfloor \frac{kn}{2k+1}\rfloor -j$ into $\cG$. Proving this for every possible $m$ would complete the proof of the upper bound on $w_n(G)$.

Now for $h< n/2$, we have $\frac{w_n(h)}{w_n(h-1)}=\frac{\binom{n}{h}}{\binom{n}{h-1}}=\frac{n-h+1}{h}$. So if $h< \frac{kn}{2k+1}$, then $\frac{w_n(h)}{w_n(h-1)}>\frac{k+1}{k}$. Hence writing $z:=\lfloor \frac{kn}{2k+1}\rfloor$, we have 
\[
\frac{w_n(m(j))}{w_n(z-j)}\ge \frac{w_n(z)}{w_n(z-j)}>\left(\frac{k+1}{k}\right)^j>4k^2,
\]
if $j\ge 2(k+1)\log (2k)$ as $(1+\frac{1}{k})^{k+1}>e$ for all $k\le 1$. This means that for any $m$ such that $m(j)=m$ implies $j\ge 2(k+1)\log (2k)$, we have
\[
\sum_{G\in\cG: |G|=m}w_n(G)+\sum_{G\in \cG:|G|=z-j, m(j)=m}w_n(G)\le n\binom{n}{m}.
\]
This implies $\sum_{G\in \cG}w(G)\le n\sum_{j\ge z-2(k+1)\log (2k)}\binom{n}{j}$, as needed.

\medskip

The second statement of the theorem follows as the first statement and Proposition \ref{constr} imply $\binom{n}{\le \lfloor\frac{kn}{2k+1}\rfloor-(k+1)\log 2k}\le 2^n-\vex(n,C_{2k+1})\le \binom{n}{\le \lfloor\frac{kn}{2k+1}\rfloor}$, and for any $m\le (1/2-\varepsilon)n$ we have $\binom{n}{\le m}\le C_\varepsilon \binom{n}{m}$.
\end{proof}

\section{Statements and Declarations}

Research supported by NKFIH under grants FK132060 and KKP-133819.

The authors have no relevant financial or non-financial interests to disclose.

\end{document}